\documentclass[11pt]{amsart}

\usepackage{amssymb}
\usepackage{amsmath, amsfonts}
\usepackage{amsthm}
\usepackage{epsfig}
\usepackage{graphicx}
\usepackage{graphics}
\usepackage{float}
\usepackage{subfigure}
\usepackage{multirow}
\usepackage{color}
\usepackage{lineno}
\usepackage{fullpage}
\usepackage[normalem]{ulem} 
\usepackage{xspace}
\usepackage{wrapfig}
\usepackage{amsmath,amsfonts,amssymb,amscd,amsthm,amsbsy,epsf}
\usepackage{euscript}
\usepackage{url}
\usepackage[colorlinks=true,citecolor=blue]{hyperref}

\newtheorem{theorem}{Theorem}
\newtheorem{meta-thm}[theorem]{Meta-Theorem}
\newtheorem{lemma}[theorem]{Lemma}

\newtheorem{remark}[theorem]{Remark}

\newtheorem{conjecture}[theorem]{Conjecture}

\newcommand\elem{ \end{lemma} }
\newcommand\lem[1]{ \begin{lemma}\label{#1}}

\newcommand\beq[1] { \begin{equation}\label{#1} }
\newcommand{\eeq}{ \end{equation} }

\newcommand\beqa[1]{ \begin{eqnarray} \label{#1}}
\newcommand{\eeqa}{ \end{eqnarray} }
\newcommand{\beqano}{ \begin{eqnarray*} }
\newcommand{\eeqano}{ \end{eqnarray*} }

\newcommand{\balino}{ \begin{align*} }
\newcommand{\ealino}{ \end{align*} }

\definecolor{indiagreen}{rgb}{0.07, 0.53, 0.03}

\def\e{\varepsilon}

\def\D{{\mathcal D}}

\def\M{{\mathcal M}}

\def\O{{\mathcal O}}

\def\integer{{\mathbb Z}}

\def\real{{\mathbb R}}

\def\r{{\rho}}

\def\t{{\theta}}

\newcommand{\Addresses}{{
\bigskip
\footnotesize

A.P.~Bustamante, \textsc{School of Mathematics, Georgia Institute of Technology}\par\nopagebreak \textit{E-mail address}, \texttt{apb7@math.gatech.edu}

\medskip

R.C.~Calleja, \textsc{Department of Mathematics and Mechanics IIMAS, National Autonomous University of Mexico (UNAM)}\par\nopagebreak \textit{E-mail address}, \texttt{calleja@mym.iimas.unam.mx}

}}

\begin{document}

\title{Corrigendum and addendum to \emph{"Computation of domains of analyticity for the dissipative standard map in the limit of small dissipation"}}

\author{ Adri\'an P. Bustamante \and  Renato C. Calleja } 



\begin{abstract}
We correct some tables and figures in [A.P. Bustamante and R.C. Calleja, \emph{Physica D: Nonlinear Phenomena}, 395 (2019), pp. 15-23]. We also report on the new computations that verify the accuracy of the data and extend the results. The new computations have led us to find new patterns in the data that were not noticed before. We formulate some more precise conjectures. 
\end{abstract}

\maketitle


\section{Introduction}
The goal of this note is to present a correction of some of the  tables and figures presented in  \cite{Bus-Cal-19}, see Section~\ref{correction}. We have also revised and extended the results with a new implementation of the algorithms. This has lead to some new patterns in the data (Section \ref{patterns}) and new verifications (Section \ref{validation}).

We recall that the aim of \cite{Bus-Cal-19} was to study quantitatively the domains of analyticity of quasiperiodic orbits for the dissipative standard map \eqref{dis-map} through a careful analysis of their Lindstedt series, as well as with non-perturbative computations. The results in \cite{Bus-Cal-19} agreed with the conjectures in \cite{Cal-Cel-Lla-16}. In particular, the result in \cite{Bus-Cal-19} verifies numerically the conjecture about the optimality of the domains of analyticity described in \cite{Cal-Cel-Lla-16}. The qualitatively conjectured optimal domain of analyticity for the map \eqref{dis-map} does not contain any ball with center at the origin nor angular sectors with width larger than $\pi /3$, so one does not expect the Lindstedt series to converge. 
The shape of the domain of analyticity suggest the Lindstedt expansions might  belong to a Gevrey class.  

In this work we present corrections and extended results related to the Gevrey character of the Lindstedt series that was also studied in \cite{Bus-Cal-19}. In particular, some of the figures and tables presented in \cite{Bus-Cal-19} are not accurate, see Section \ref{correction}. The corrected tables included here contain sharper results. With the corrected data  and the extended computations performed, we have reformulated a conjecture about the Gevrey character of the Lindsted series, see Conjecture \ref{conjecture1}. 

We note that some rigorous studies compatible with the conjectures in \cite{Bus-Cal-19} and Section~\ref{patterns} have been obtained recently in \cite{BustamanteL20}.

\section{Summary of \cite{Bus-Cal-19}}

\subsection{Lindstedt series} We recall that one of the goals of 
\cite{Bus-Cal-19}  was to study some properties of  the Lindstedt series of quasi-periodic orbits for the dissipative standard map $f_\e(x_n,y_n) = (x_{n+1}, y_{n+1})$
\begin{align}
x_{n+1} &= x_n + y_{n+1} \label{dis-map} \\  
y_{n+1} &= b_\e y_n +c_\e + \e V'(x_n) \nonumber
\end{align}
defined on the cylinder $\M = \mathbb{S}^1\times \real$; $b_\e = 1-\e^3$, $V'(x) = \frac{1}{2\pi} \sin(2\pi x)$. When one chooses the parameter $c_\e$ appropriately, it is known that \eqref{dis-map} has an analytic invariant circle corresponding to a quasi periodic orbit with Diophantine frequency $\omega$. 

It is known that quasi periodic orbits of  \eqref{dis-map}  can be described by a 1-periodic function $u_\e: \mathbb{S}^1\rightarrow \real $ and a constant $c_\e$ satisfying \begin{equation}\label{inv-eq}
    E_{c_\e}[u_\e ] = 0
\end{equation}
where $ E_{c_\e}[u_\e(\t)] \equiv u_\e(\t +\omega) - (1+b_\e)u_\e(\t) +b_\e u_\e(\t -\omega) +(1-b_\e)\omega -c_\e +\e V'(\t +u_\e(\t)) $.

Lindstedt expansions are obtained by considering the formal power series $u_\e(\t) = \sum_{k=0}^\infty u_k(\t)\e^k $ and $c_\e = \sum_{k=0}^\infty  c_k\e^k$, and solving \eqref{inv-eq} order by  order. 
The coefficients $u_k$  and  $c_k$ are determined by the following cohomology equation \begin{equation} \label{order-equa}
    L_\omega u_k(\t) -c_k +u_{k-3}(\t) -u_{k-3}(\t -\omega) = S_k(\t), \qquad k\geq 4
\end{equation} 
where $L_w \varphi(\t) \equiv \varphi(\t+\omega) -2\varphi(\t) +\varphi(\t - \omega) $, and  $\e V' \equiv \sum_{k=0}^\infty S_k(\t)\e^k$.
\\

We note that given that the Lindstedt series in this case are not convergent in any ball and the terms grow very fast, the numerical calculation of the coefficients $u_k$ is much more unstable that in the cases where the Lindstedt series converges.

\subsection{Gevrey character of Linsdtedt series}\label{sec-gevrey}
One of the goals of \cite{Bus-Cal-19} was to study, numerically, the Gevrey character of the Lindstedt series $\sum u_k\e^k$. To do this we considered the quantities \begin{equation}\label{growths} 
A_\r(k) \equiv \frac{1}{k}\log \| u_k \|_{\r}, \qquad  H^r(k) \equiv \frac{1}{k}\log \| u_k\|_{W^r}
\end{equation}
which measure the growth of the coefficients of the Lindstedt series using different norms. The norms were chosen as  $\|f \|_{\rho} = \sum_{\ell \in \integer}|\hat{f}_\ell |^2 e^{2\pi|\ell|\rho} $ and $ \|f\|_{W^r}^2 = \sum_{k\in \integer} (2\pi k)^{2r}|f_k|^2 $.

We recall that a formal power series, $\sum f_n\e^n$, is $\sigma$-Gevrey with respect to a norm, $\|\cdot \| $, if the coefficients satisfy $$\|f_n\|\leq CR^n n^{\sigma n}.$$ Equivalently, $$\frac{1}{n}\log \|f_n\|  \sim \sigma \log(n) + \log(R) $$ for $n$ large enough.

\section{Correction to \cite{Bus-Cal-19} } \label{correction}

The main correction is that the 
data in Table 1, Table 2 and in the plots on Figure 2, and Figure 3 in \cite{Bus-Cal-19} do not correspond to their labels. These tables and figures were included to study the growth of the coefficients of the Linsdtedt series $\sum u_k \e^k$, corresponding to the frequency $\omega = \frac{\sqrt{5} -1 }{2}$.

The correct Table 1 and Figure 2 in \cite{Bus-Cal-19}  must be:

\begin{table}[H]
\begin{tabular}{ |p{3cm}||p{2cm}|p{2cm}|p{2cm}| }
 \hline
 \multicolumn{4}{|c|}{ $e_{\rho}(k): = A_\rho (k) -( \log(R) + \sigma \log(k))$  } \\
 \hline
 &$ R$  & $\sigma$ & $\|e_{\rho}\|_\infty$  \\
 \hline
$\rho = 0.1$ &  0.672269  & 0.227899 & 0.020793  \\
$\rho = 0.01$ & 0.585740 & 0.238324 & 0.019491\\
$\rho = 0.001$ & 0.576278 & 0.240049  & 0.019325 \\
 $\rho = 0.0001$ &  0.575333 & 0.240225 & 0.019280 \\
 $\rho = 0.00001$ &   0.575239  & 0.240243 & 0.019282 \\
 $\rho = 0.000001$ &  0.575230 &  0.240244 & 0.019279 \\
 $\rho = 0.0000001$ & 0.575229 & 0.240244 & 0.019278 \\
 
 \hline
\end{tabular}
\caption{Numerical fit of a function $\log(R) +\sigma\log(k)$ to the data  $A_\rho(k)$ for different values of $\rho$ and  frequency $\omega = \frac{\sqrt{5}-1}{2}$. Computations were done using $2^{13}$ Fourier coefficients and 600 digits of precision. The numerical fit was made in for $100\leq k\leq 300$.}

\label{tab-analytic}
\end{table}

\begin{figure}[H]
\includegraphics[ trim = 8cm 10cm 8cm 10cm, width=4truecm]{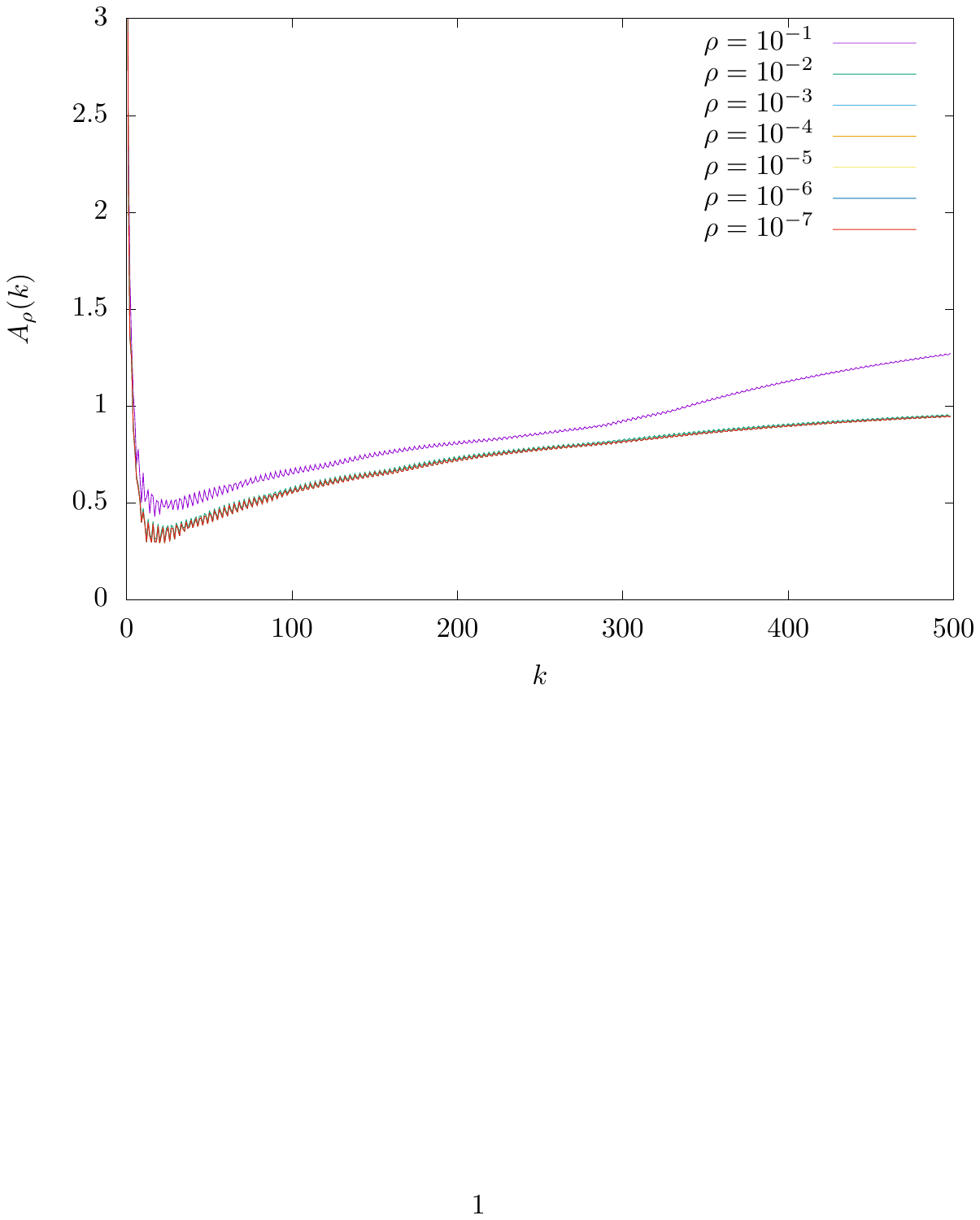}
\caption{ Plot of $A_\r(k)$, $1\leq k\leq 500$, for the frequency $\omega = \frac{\sqrt{5}-1}{2}$.  } 
\label{norm1}
\end{figure}

The corrected Table 2 and Figure 3 in \cite{Bus-Cal-19} must be:

\begin{table}[H]
\begin{tabular}{ |p{2cm}|p{2cm}|p{2cm}|p{2cm}| }
 \hline
 \multicolumn{4}{|c|}{ $e_{r}(k) := H^r (k) -( \log(R) + \sigma \log(k))$  } \\
 \hline
 &$ R$ & $\sigma$ & $\|e_{r} \|_\infty $   \\
 \hline
 $r = 1$ & 0.685071 & 0.212840  & 0.020144 \\
 $r = 2$ & 0.816610 & 0.185284 & 0.023905 \\
 $r = 3$ & 0.974288 & 0.157572 & 0.028145 \\
 $r = 4$ & 1.163403 & 0.129713 & 0.032216  \\
 $r = 5$ & 1.390238 & 0.101731 & 0.036129  \\
 $r = 6$ & 1.662287 & 0.073651 & 0.039905 \\
 
 \hline
\end{tabular}
\caption{Numerical fit of a function $\log(R) +\sigma\log(k)$ to the data $H^r(k)$ for different values of $r$ and frequency $\omega = \frac{\sqrt{5}-1}{2}$. Computations were done using $2^{13}$ Fourier coefficients and 600 digits of precision. The numerical fit was made for $100\leq k\leq 300$.}

\label{tab-analytic2}
\end{table}
 
 \begin{figure}[H]
\begin{center}
\includegraphics[ trim = 8cm 10cm 8cm 10cm, width=4truecm]{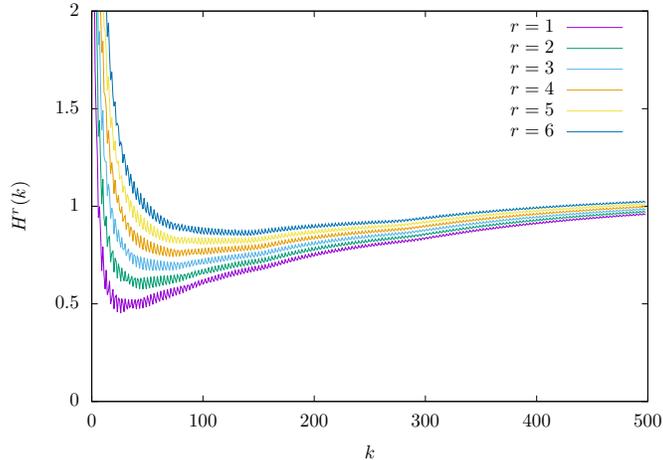}
\caption{ Plot of $H^r(k)$, $1\leq k\leq 500$, for the frequency $\omega = \frac{\sqrt{5}-1}{2}$.  } 
\end{center}
\end{figure}

We note that the numbers $R$ and $\sigma$ in Table \ref{tab-analytic} and Tale \ref{tab-analytic2} are just the raw numbers obtained by fitting numerically functions of the form $\log(R) + \sigma \log(k)$ to the data $A_\r(k)$ and $H^r(k)$, we are not sure how to assess the reliability of these numbers. Also, we have added a column with a measure of the remainder,  $\|e\|_\infty$, between the numerical fit and the data, this column was not included in the tables in \cite{Bus-Cal-19}. The measure of these remainders, which looks a little bit worrisome, seems to come from an \emph{oscillatory} behavior in the data, the structure of the remainders is studied in Section \ref{new-patterns}.

The problem with the figures and tables in \cite{Bus-Cal-19} is that the data used on them corresponded to the quantities $\frac{1}{k}\log\|k!u_k\|$ and not to the data given by $A_\rho(k)$ and $H^r(k)$, defined in \eqref{growths}. Note that, by the well known Stirling's formula $\log(k!) = k\log(k) -k +O(\log(k)) $, if $k^{-1}\log\|u_k\|\approx \log(R) + \sigma\log(k) $ then $k^{-1}\log\|k!u_k\| \approx \log(\Tilde{a}) + (\sigma +1)\log(k)$ for $k \gg 1$.  We recall that the values corresponding to the column $\sigma$ on the tables in \cite{Bus-Cal-19} gave $\sigma  \approx 1$. The fact that the values of $\sigma$ in Table \ref{tab-analytic} and Table \ref{tab-analytic2} are not approximate to zero can be explained by how the numerical fits are done this time, which is explained in the next paragraph.

\begin{figure}[H]

\includegraphics[ trim = 8cm 10cm 8cm 10cm, width=4truecm]{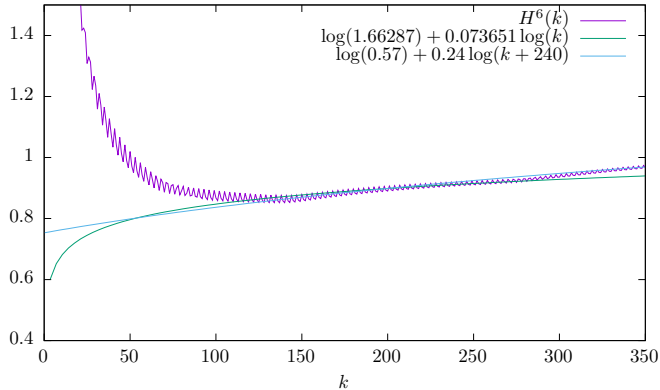}
\caption{ Comparison between $H^6(k)$ and  two different numerical fits for  $100\leq k \leq 300$, $\omega = \frac{\sqrt{5}-1}{2}$. It can be observed that introducing  a translation, $b$, could make a significant change in the exponent $\sigma$. } 
\label{logs}
\end{figure}

We note first that \cite{Bus-Cal-19} used fits of the form $\log(R) + \sigma \log(k+b)$ which involve an extra parameter $b$. We consider that omitting the translation by $b$ is more suitable for a systematically study of the growth of the coefficients of the Lindstedt series, see Section \ref{sec-gevrey}. Note that by adjusting $R$ and $b$ one can get $\log(\tilde{R}) + \tilde{\sigma} \log(k) \approx \log(R) + \sigma\log(k+b)$ with $k_1\leq k \leq k_2$, for example, $\log(1.66287) + 0.073651\log(k) \approx \log(0.57) + 0.24\log(k+240)$ for $100\leq k\leq 300$, see Figure \ref{logs}.  We also note the numerical fits are  made  taking a smaller range for $k$ in $A_\rho(k)$ and $H^r(k)$, the range for $k$ in \cite{Bus-Cal-19} was $100\leq k\leq 1000$ (which is another reason for which the tables needed to be corrected). The fits in Table \ref{tab-analytic} and Table \ref{tab-analytic2} were made considering $100\leq k \leq 300$, this is due to the fact that we consider that the errors, in the computation of the coefficients $u_k$, are small enough within this range of parameters, see Section \ref{validation}. Finally, the factor $O(\log(k)/k)$, that one gets using  Stirling's formula, satisfies $O(\log(k)/k)=O(10^{-2}) $ for $100\leq k \leq 300$, which could also affect the values of $\sigma$ at  order $10^{-1}$. We consider that the observations above explain why the values in the column $\sigma$ in Table \ref{tab-analytic} and Table \ref{tab-analytic2} are not only translations, by 1, of the values  obtained in \cite{Bus-Cal-19}.

\begin{figure}[H]

\includegraphics[ trim = 8cm 10cm 8cm 10cm, width=4truecm]{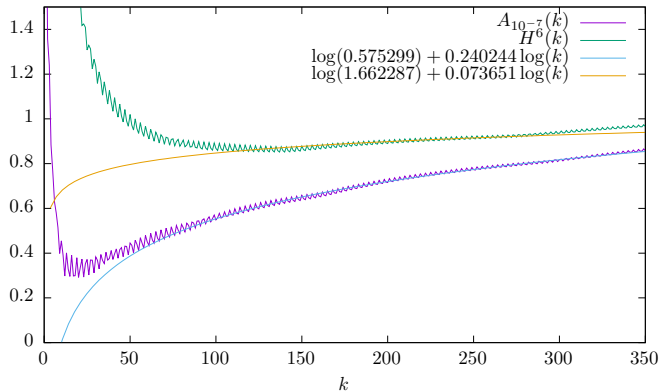}
\caption{ Comparison between $H^6(k)$ and $ A_{10^{-7}}(k)$ with their respective numerical fits, $\omega = \frac{\sqrt{5}-1}{2}$.  } 
\label{different-fits}
\end{figure}

For the sake of completeness we include a comparison between $A_{10^{-7}}(k)$, $H^6(k)$ and their respective numerical fits, see Figure \ref{different-fits}. Note that even if the norms considered in $A_\r(k)$ and $H^r(k)$ are in principle not compatible, the fact that $A_\r(k)$ and $H^r(k)$  have similar trends seems to indicate that there is a mechanism which is captured for any norm for the functions we study. This suggest that a more detailed study of the structure of this functions could be interesting.

The corrections in the tables have important consequences for the statement of Conjecture 9 in \cite{Bus-Cal-19}. The conjecture proposed the Gevrey character of the Lindstedt expansions with Gevrey exponent $\sigma \approx 1$, according to the tables included in \cite{Bus-Cal-19}. The corrected tables, Table \ref{tab-analytic} and Table \ref{tab-analytic2} (with more reliable data), suggest that the conjecture about the Gevrey character is still true but with a different exponent $\sigma$. We reformulate the conjecture in Section \ref{more-frequencies}, after we present the results we have obtained with the extended computations that have been performed. We recall that the computations in \cite{Bus-Cal-19} were done only for the frequency $\omega = \frac{\sqrt{5}-1 }{2} $, in the next section we present also results for different values of $\omega$.

\section{Some new patterns and extension of the computations}\label{patterns}
Since the publication of \cite{Bus-Cal-19} we have run several modifications of the program and re-implemented some of the algorithms. 
This allowed us to find some new patterns in the data and extend
the computations to other frequencies $\omega$. The new results allow us to reformulate the conjecture established in \cite{Bus-Cal-19}, Conjecture \ref{conjecture1}, and also give evidence of new patterns that were not noticed before, see Conjecture \ref{conjecture2}.



\subsection{Results for different frequencies} \label{more-frequencies}

We recall that the computations in \cite{Bus-Cal-19} were done using the Frequency $\omega = \frac{\sqrt{5}-1}{2}$. This time, we have performed the computations also for different frequencies, of the same Diophantine type, and we have found a \emph{similar} behavior in the growth of  the coefficients of the Lindstedt series. We present the results below.

\begin{figure}[H]
\includegraphics[ trim = 7cm 10cm 7cm 10cm, width=6truecm ]{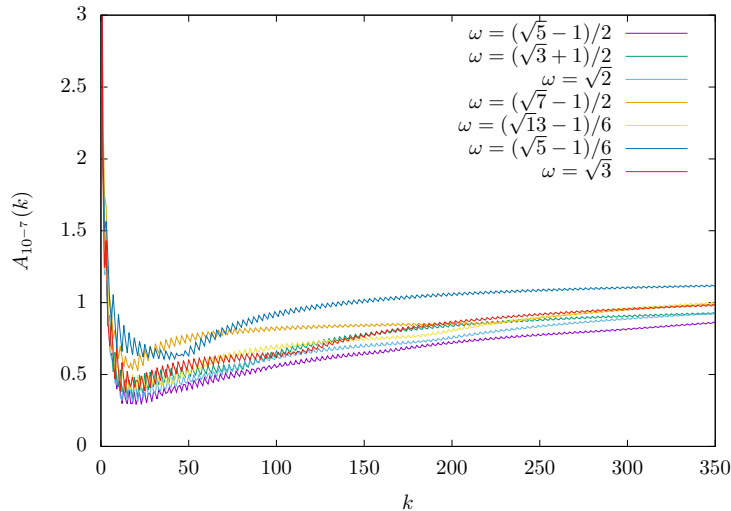}
\caption{Graph of $A_\rho(k)$ for different values of the frequencies $\omega$, $\r = 10^{-7}$. }
\label{all-freq}
\end{figure}

Figure \ref{all-freq} contains a plot of $A_\rho(k)$, $\rho = 10^{-7}$, for all the frequencies considered.

The plots in Figure \ref{all-freq} seems to indicate a logarithm growth for all the frequencies considered. To study more systematically the growth of the coefficients of the Lindsted series we have also fitted numerically functions of the form $\log(a) + \sigma\log(k)$, the results are summarized in Table \ref{tab-fits}. We note that all the frequencies considered belong to the same Diophatine class $\D(\nu,1)$, where $\omega\in\D(\nu,\tau)$ means that $|e^{2\pi i k\omega} -1|\geq v|k|^{-\tau} $.

\begin{table}[H]
\begin{tabular}{ |p{6cm}||p{2cm}|p{2cm}|p{2cm}| }
 \hline
 \multicolumn{4}{|c|}{$e_{\omega}(k) = A_\rho (k) - \log(R) + \sigma \log(k),\qquad \rho = 10^{-7}$ } \\
 \hline
 & $R$ & $\sigma$ & $\|e_{\omega} \|_\infty $ \\
 \hline
 $ \omega = \frac{\sqrt{5}-1}{2} = [0,1,1,1,1,1,1,...] $ &  0.575229 & 0.240244 & 0.019278 \\
 $ \omega = \frac{\sqrt{3}-1}{2} = [0,2,1,2,1,2,1,...]  $ &   0.695887  & 0.225349 & 0.047762 \\
 $\omega = \sqrt{2} = [1,2,2,2,2,2,2,...] $ &  0.583365 &  0.247799 & 0.033104 \\
 $\omega = \sqrt{3} = [1,1,2,1,2,1,2,1,...] $ & 0.460186  &  0.307029 & 0.038801 \\
 $ \omega = \frac{\sqrt{7}-1}{2} = [0,1,4,1,1,4,1,1,...]  $ &   1.300597  & 0.112924 & 0.045704 \\
 $ \omega = \frac{\sqrt{13}-1}{6} = [0,2,3,3,3,3,3,...]  $ &   0.582937  & 0.258504 & 0.047840 \\
 $ \omega = \frac{\sqrt{5}-1}{6} = [0,4,1,5,1,5,1,5,...]  $ &   1.235768  & 0.158503 & 0.042327 \\
 \hline
\end{tabular}

\caption{Numerical fit of a function $\log(R) +\sigma\log(k)$ to the data  $A_\rho(k)$ for different values of the frequency $\omega$ and $\rho = 10^{-7}$. Computations were done using $2^{13}$ Fourier coefficients and 600 digits of precision. The numerical fit was made in for $100\leq k\leq 300$.}
\label{tab-fits}
\end{table}

\begin{figure}[H]
\includegraphics[trim = 8cm 10cm 8cm 10cm,width=3.5truecm]{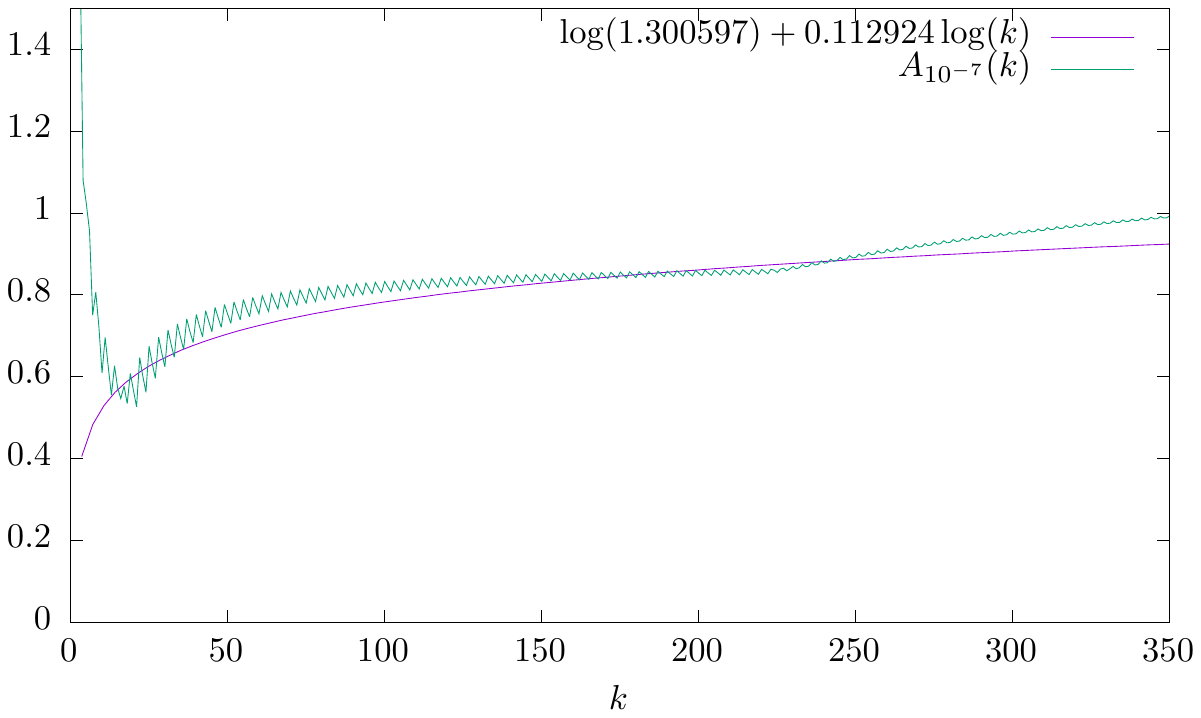}\hfil\hfil\hfil
\includegraphics[trim = 8cm 10cm 8cm 10cm,width=3.5truecm]{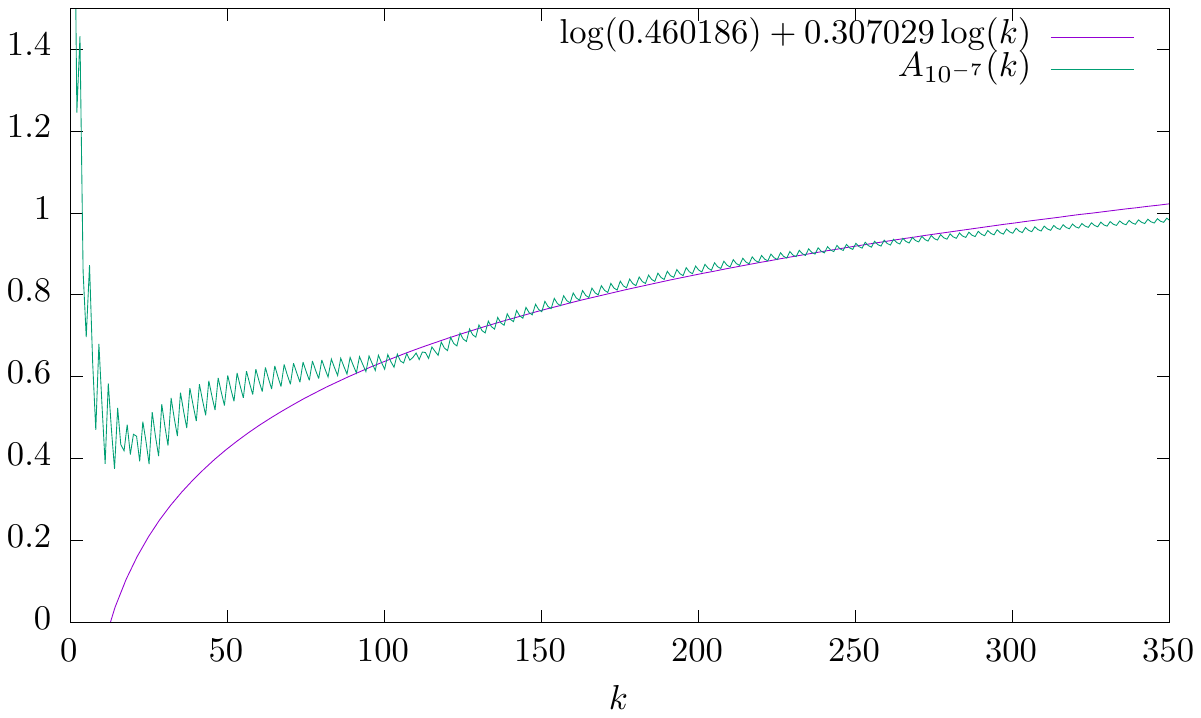}
\caption{Comparison between  $A_\r(k)$, $\r= 10^{-7}$, and its corresponding numerical fit . Left panel: values for the frequency $\omega = \frac{\sqrt{7}-1}{2}$. Right panel: values for the frequency $\omega=\sqrt{3}$ }
\label{norm_other_omegas}
\end{figure}

\begin{figure}[H]
\includegraphics[trim = 8cm 10cm 8cm 10cm,width=3.5truecm]{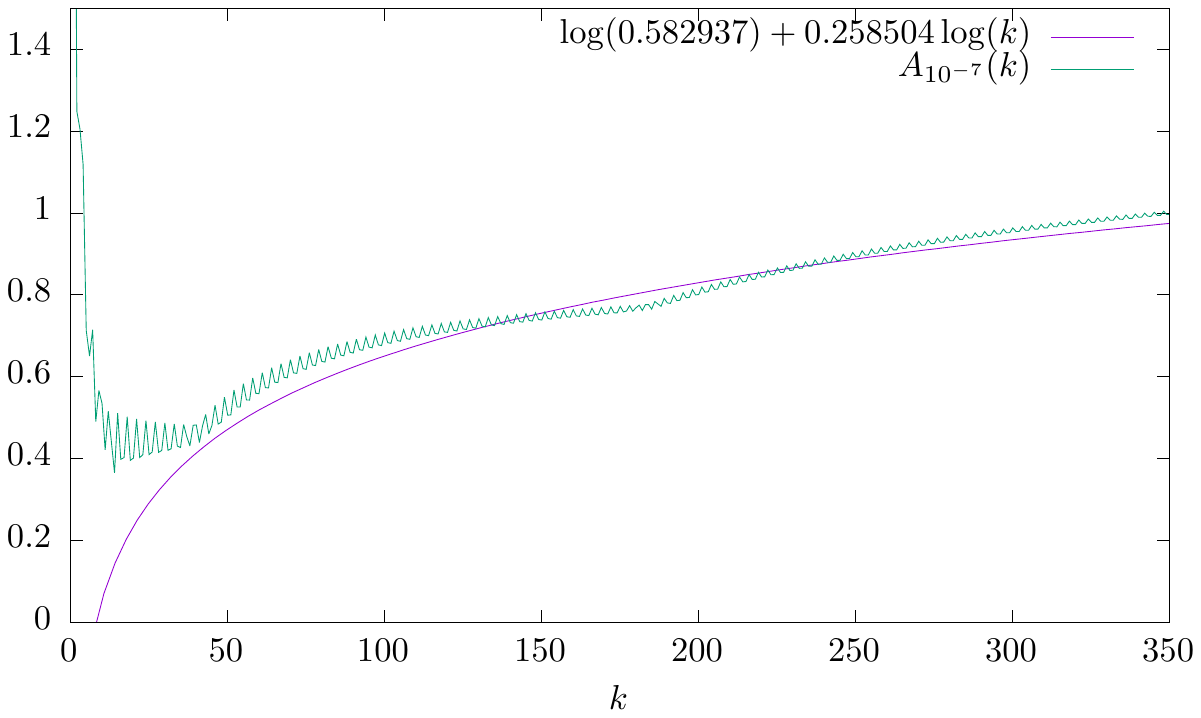}\hfil\hfil\hfil
\includegraphics[trim = 8cm 10cm 8cm 10cm,width=3.5truecm]{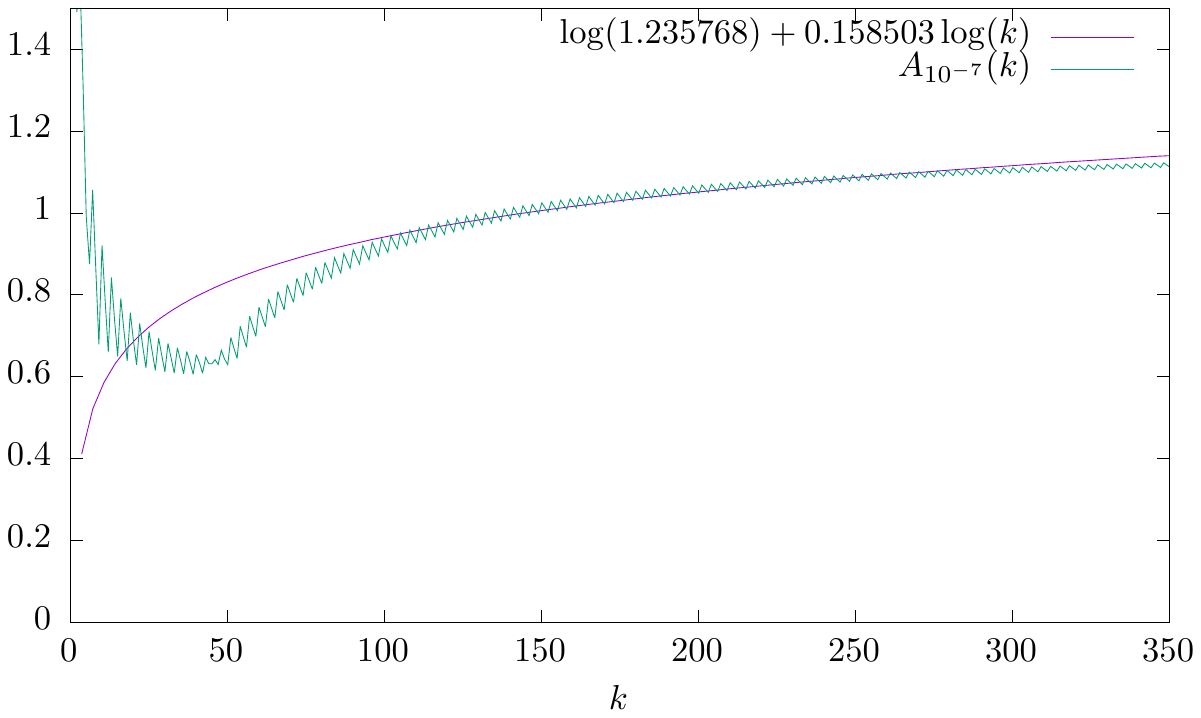}
\caption{ Comparison between  $A_\r(k)$, $\r= 10^{-7}$, and its corresponding numerical fit. Left panel: values for the frequency $\omega = \frac{\sqrt{13}-1}{6}$. Right panel: values for the frequency $\omega= \frac{\sqrt{5}-1}{6}$ }
\label{norm_other_omegas2}
\end{figure}

Figure \ref{norm_other_omegas} and Figure \ref{norm_other_omegas2} contain comparisons between the quantities $A_\r$ and their respective numerical fits.

The extension of the computations to other frequencies and the information summarized in Table \ref{tab-fits} allow us to reformulate Conjecture 9 in \cite{Bus-Cal-19} for a more general case.

\begin{conjecture} \label{conjecture1}
Given $\omega\in\D(\nu,1)$, the Lindstedt series, $u_\e =\sum u_k\e^k$, of quasi-periodic orbits for the map \eqref{dis-map}  belongs to a Gevrey class with Gevrey exponent $\sigma \leq 0.307$. That is, $\|u_n\|_\rho \leq CR^n n^{\sigma n}$ with $\sigma \leq 0.307$ and $\rho\leq 10^{-7}$.  
\end{conjecture}

\begin{remark}
It is worth to note that Conjecture \ref{conjecture1} is compatible with the rigorous results obtained in \cite{BustamanteL20}. Considering the map \eqref{dis-map}, with dissipation $b_\e = 1-\e^3$ and a frequency $\omega\in \D(\nu,1)$, the rigorous results in \cite{BustamanteL20} yield a Gevrey exponent $\sigma = 2/3$.

It is also important to note that the results in \cite{BustamanteL20} give the same upper bound of the  Gevrey exponent for frequencies, $\omega$, of the same Diophantine type $\D(\nu,\tau)$. The behavior observed in Figure \ref{all-freq}, Figure \ref{norm_other_omegas}, Figure \ref{norm_other_omegas2}, and Table \ref{tab-fits} seems to indicate that the upper bound of the Gevrey exponent found in \cite{BustamanteL20} is not optimal, but seems to be  within a factor 2 for being optimal. 
%
\end{remark}

\subsection{New patterns} \label{new-patterns}
\begin{figure}[H]
\includegraphics[ trim = 7cm 10cm 7cm 10cm, width=6truecm ]{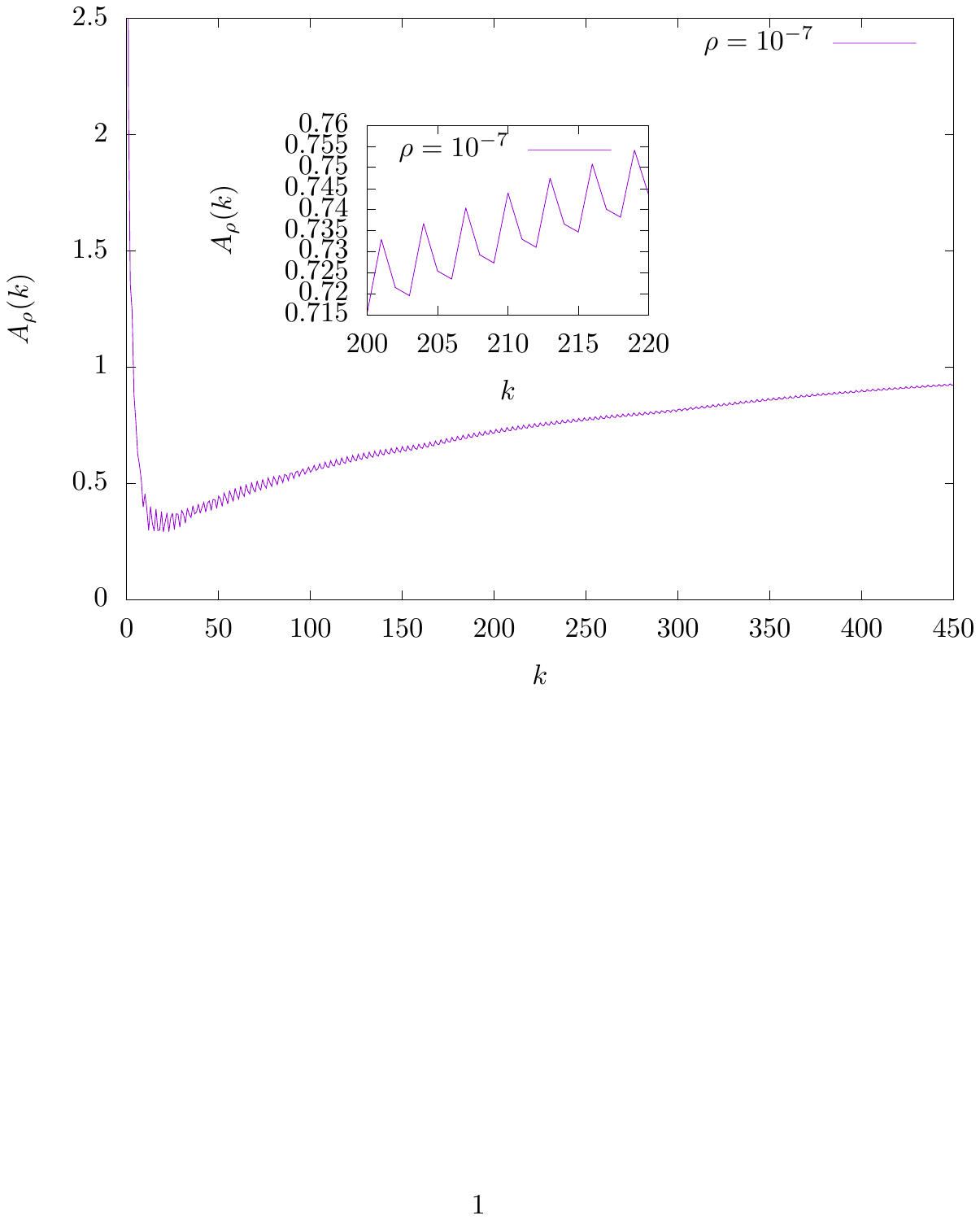}
\caption{Graph of $A_\rho(k)$ for $\r = 10^{-7}$, $\omega = \frac{\sqrt{5}-1}{2}$. The same oscillatory behavior is also present for $H^r(k)$. }
\label{zoom}
\end{figure}

A careful inspection of Figure \ref{norm1}, Figure \ref{all-freq}, Figure \ref{norm_other_omegas}, and Figure \ref{norm_other_omegas2} shows that the graphs of $A_\rho(k)$ present an \emph{oscillatory behavior} of period three, see Figure \ref{zoom}. These \emph{oscillations} are present for all the values of the frequencies considered.

As we mentioned before, the coefficients of the Lindstedt series are determined by solving equation \eqref{order-equa} in which the coefficient of order $k$ depends explicitly on the coefficient of order $k-3$. This is due to the power three of $\varepsilon$ in the function $b_\varepsilon$. At the same time this phenomenon is independent of the frequency $\omega$ we choose. This gives an explanation of the appearance of an oscillating pattern observed in the inset of Figure \ref{zoom} which appears for all the frequencies we considered. However, the computations show that the amplitude of the oscillations decreases as $k$ grows and this oscillating effect fades away.

 
 To study how the amplitude of the oscillations decreases we have \emph{centralized} the oscillations by considering the differences between $A_\rho(k)$ and some moving averages. More precisely, denoting $a_k = A_\rho(k)$, $\rho = 10^{-7}$, we have considered the following \emph{centralizations} \begin{equation}
     x_k = a_k -\frac{1}{5} \sum_{j=k-2}^{k+2}a_j, 
     \qquad z_k = a_k -\frac{1}{3k} \sum_{j=k}^{k+2}ja_j, 
 \end{equation}
 
\begin{figure}[H]
\includegraphics[trim = 8cm 10cm 8cm 10cm, width=3.8truecm]{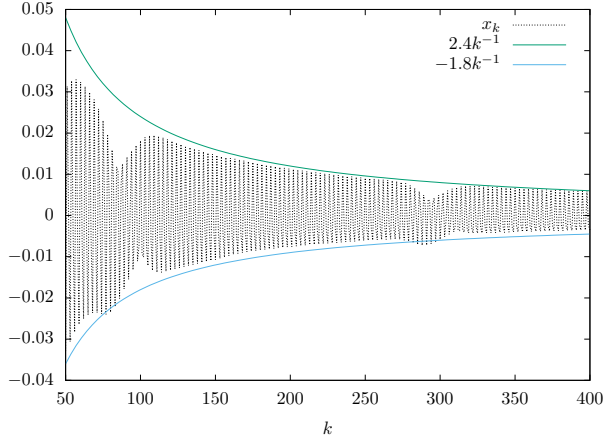}\hfil\hfil\hfil
\includegraphics[trim = 8cm 10cm 8cm 10cm, width=3.8truecm]{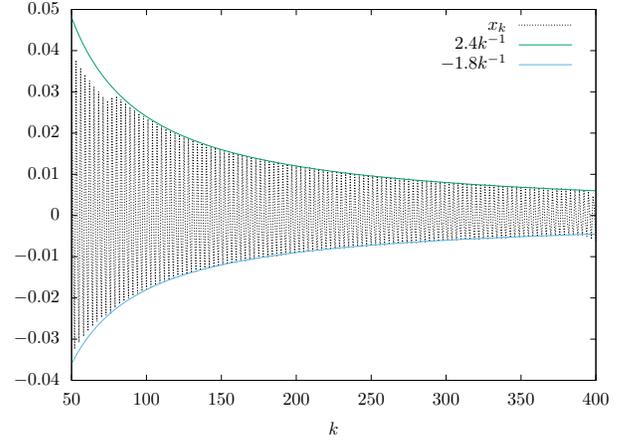}

\caption{Plots of the centralization $x_k$. Left panel: Plot for the frequency  $\frac{\sqrt{5}-1}{2}$. Right panel: Plot for frequency $\frac{\sqrt{3}-1}{2}$. } 
\label{centra-1}
\end{figure}

\begin{figure}[H]
\includegraphics[trim = 8cm 10cm 8cm 10cm, width=3.8truecm]{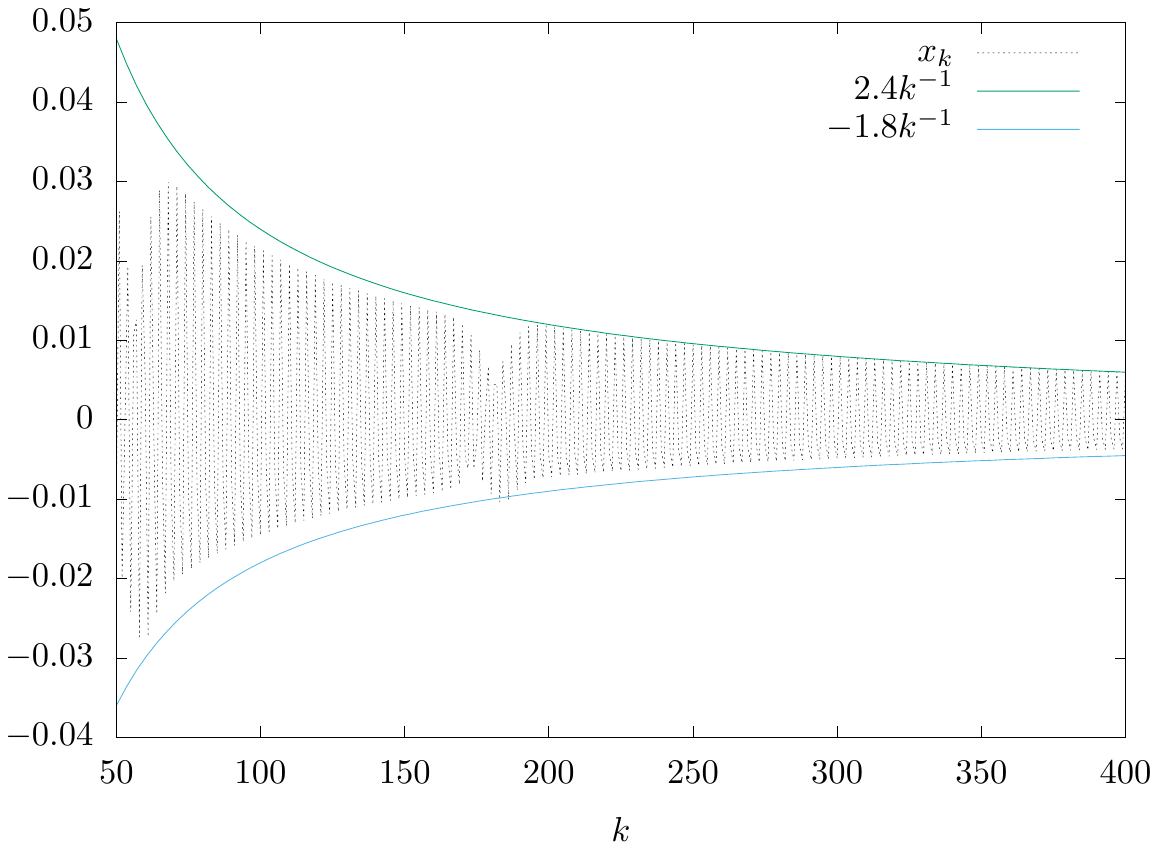}\hfil\hfil\hfil
\includegraphics[trim = 8cm 10cm 8cm 10cm, width=3.8truecm]{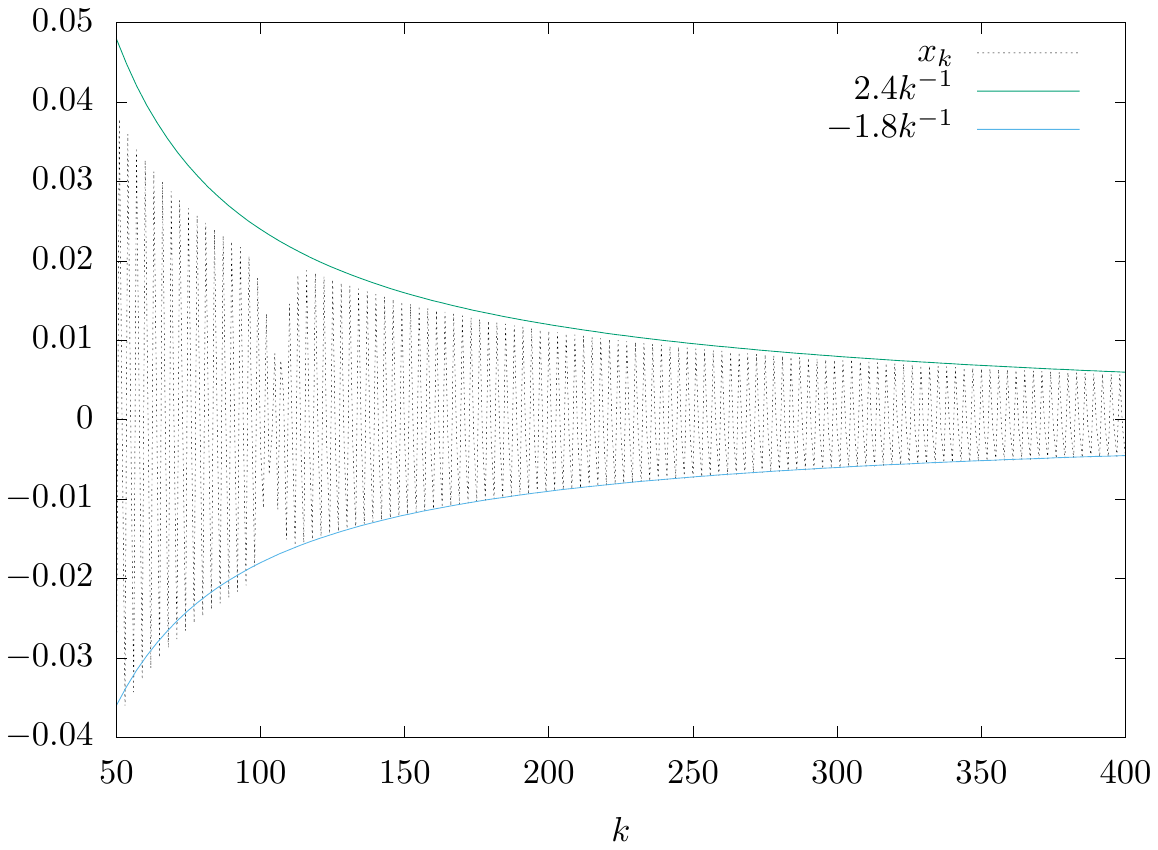}

\caption{Plots oft the centralization $x_k$. Left panel: Plot for the frequency  $\sqrt{2}$. Right panel: Plot for frequency $\sqrt{3}$. } 
\label{centra-2}
\end{figure}

\begin{figure}[H]
\includegraphics[trim = 8cm 10cm 8cm 10cm, width=3.8truecm]{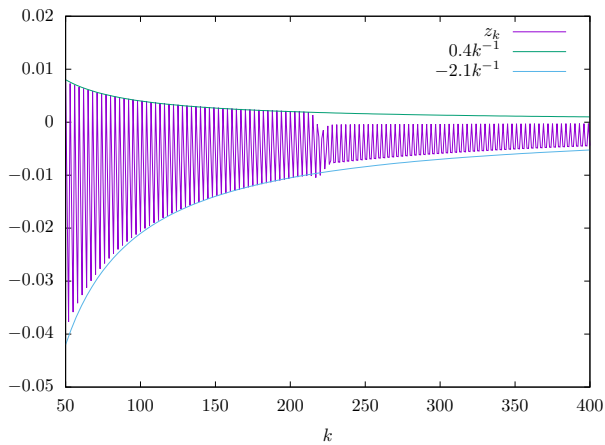}\hfil\hfil\hfil
\includegraphics[trim = 8cm 10cm 8cm 10cm, width=3.8truecm]{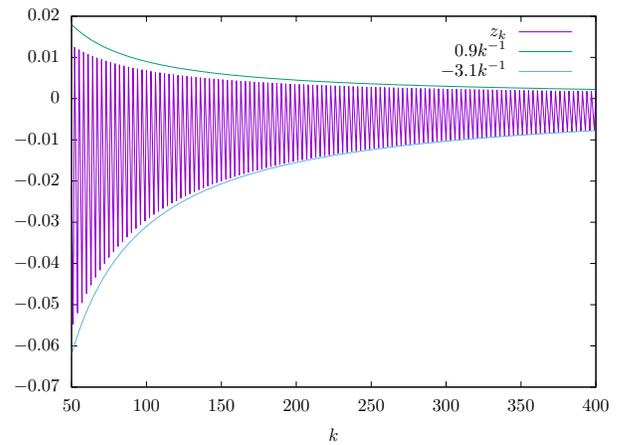}

\caption{Plots of the centralization $z_k$. Left panel: Plot for the frequency  $\frac{\sqrt{7}-1}{2}$. Right panel: Plot for frequency $\frac{\sqrt{5}-1}{6}$. } 
\label{centra-3}
\end{figure}

Since the oscillations have period three, the \emph{centralization} $x_k$ is made by subtracting a moving average that captures two periods of the oscillation. The results for $x_k$ are summarized in Figure \ref{centra-1} and Figure \ref{centra-2}.   For all the \emph{centralizations} considered it is quite surprising that the amplitude of the oscillations seems to decrease as $k^{-\beta}$, with $\beta\approx 1$. Due to this  behavior we  consider a second \emph{centralization}, $z_k$, which assumes that the oscillations decrease as $k^{-1}$. The results for $z_k$ are summarized in Figure \ref{centra-3}. 

The results collected in the figures above suggest that the \emph{centralizations} behave like $k^{-\beta}f(k)$ with $f$ a periodic function. This observation motivates the following conjecture.
\begin{conjecture}\label{conjecture2}
Let $A_\rho(k) = \frac{1}{k}\log\|u_k\|_\rho$, then $A_\rho(k) \approx \log(R) + \sigma \log(k) + k^{-\beta}f(k) $ with $\beta \approx 1$, $f(k)$ a periodic function of period 3, and $k\gg 1$.
\end{conjecture}

\section{Validation of the results} \label{validation}
To validate the results described above we verified that the 
cohomology equation \eqref{order-equa} is satisfied at every 
order with a suitable error. We also verified, as shown in \cite{Bus-Cal-19}, that the invariance equation \eqref{inv-eq} satisfies that 
$\log_{10}( \| E_{c^{\leq N}(\varepsilon)}[ u^{\leq N}_\varepsilon]\|_\infty)\sim \O((N+1)\log_{10}(\e)) $
as long as the error is above machine precision. We recall that $  E_{c^{\leq N}(\varepsilon)}[ u^{\leq N}_\varepsilon] $ means that we evaluate the operator $E$, given in $\eqref{inv-eq}$, in the finite expansions $u_\e^{\leq N} = \sum_{k=1}^N u_k\varepsilon^k$ and $c^{\leq N}(\varepsilon) = \sum_{k=0}^N c_k\varepsilon^k $.

\begin{figure}[H]
\includegraphics[trim = 8cm 10cm 8cm 10cm, width=3.8truecm]{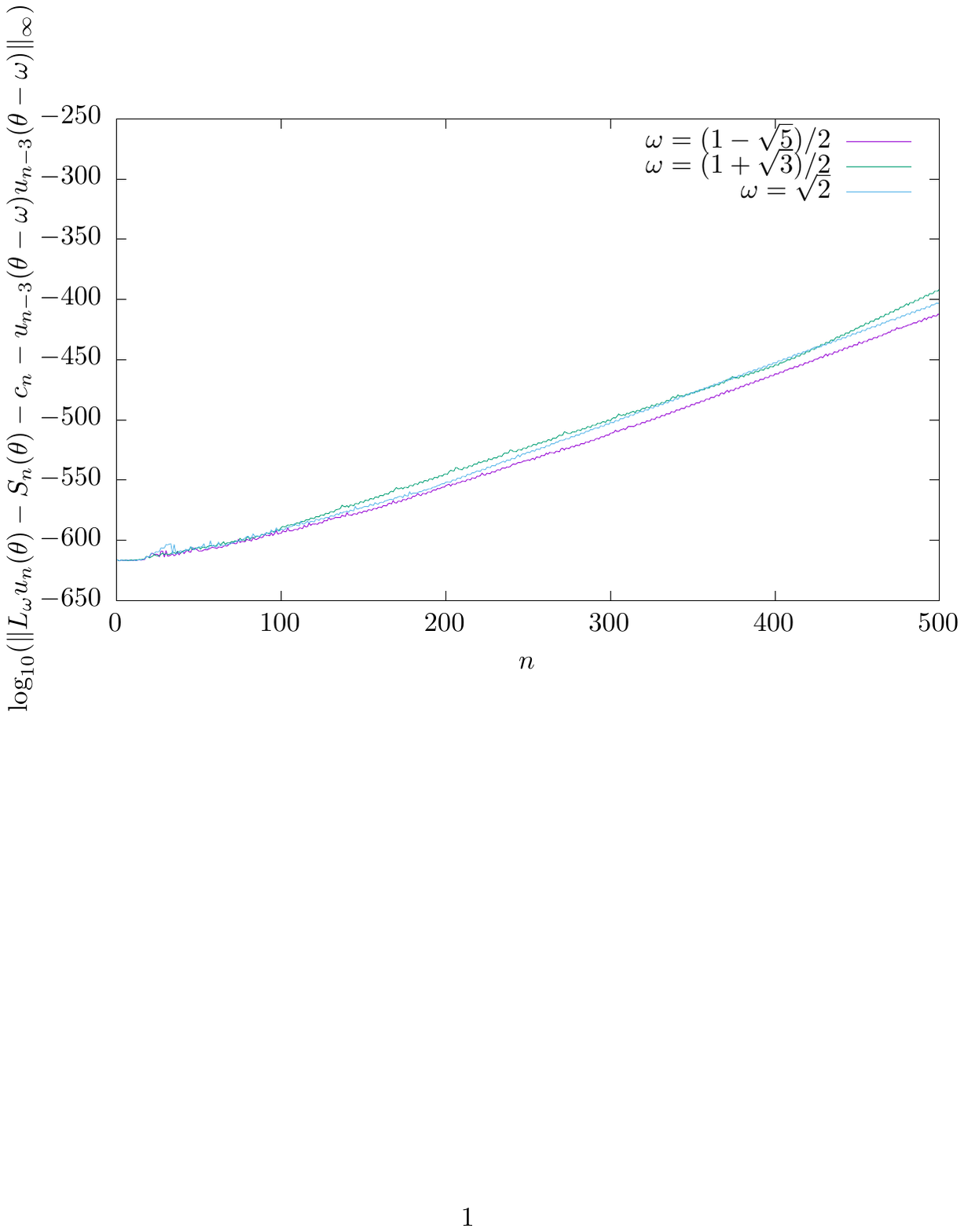}\hfil\hfil\hfil
\includegraphics[trim = 8cm 10cm 8cm 10cm, width=3.45truecm]{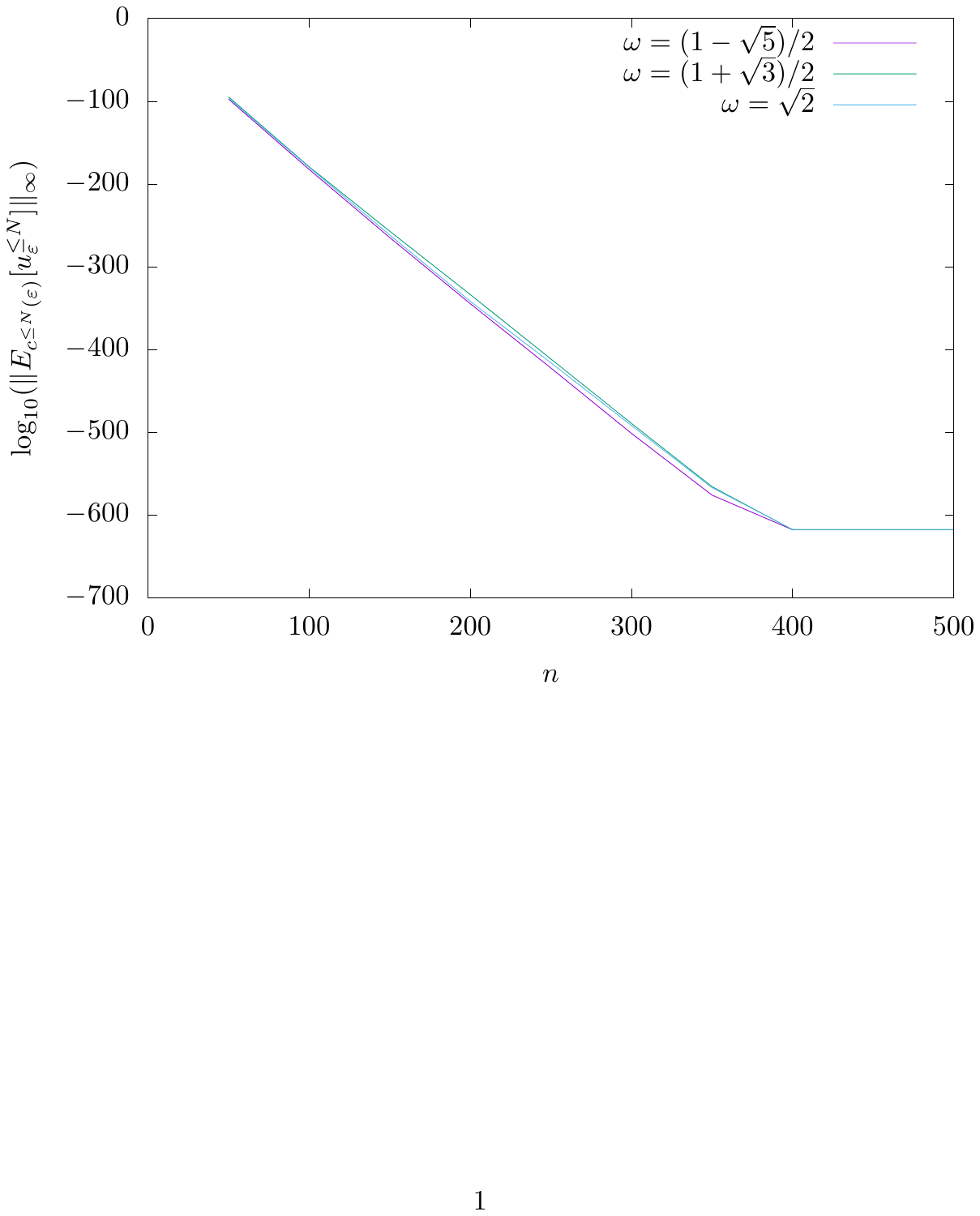}

\caption{Left panel: Plot of $\log_{10}(\| L_\omega u_n(\t) -S_n(\t) -c_n -u_{n-3}(\t -\omega) +u_{n-3}(\t) \|_\infty )$ for  different values of the frequency $\omega$, $1\leq n\leq 500$. Right panel: Plot of $\log_{10}( \| |E_{c^{\leq N}(\varepsilon)}[ u^{\leq N}_\varepsilon]\|_\infty)$, with $\epsilon =10^{-2}$.} 
\label{inveq-graphs}
\end{figure}

In Figure~\ref{inveq-graphs}, we show the results of these computations.

\begin{figure}[H]
\includegraphics[ trim = 8cm 10cm 8cm 10cm, width=4truecm]{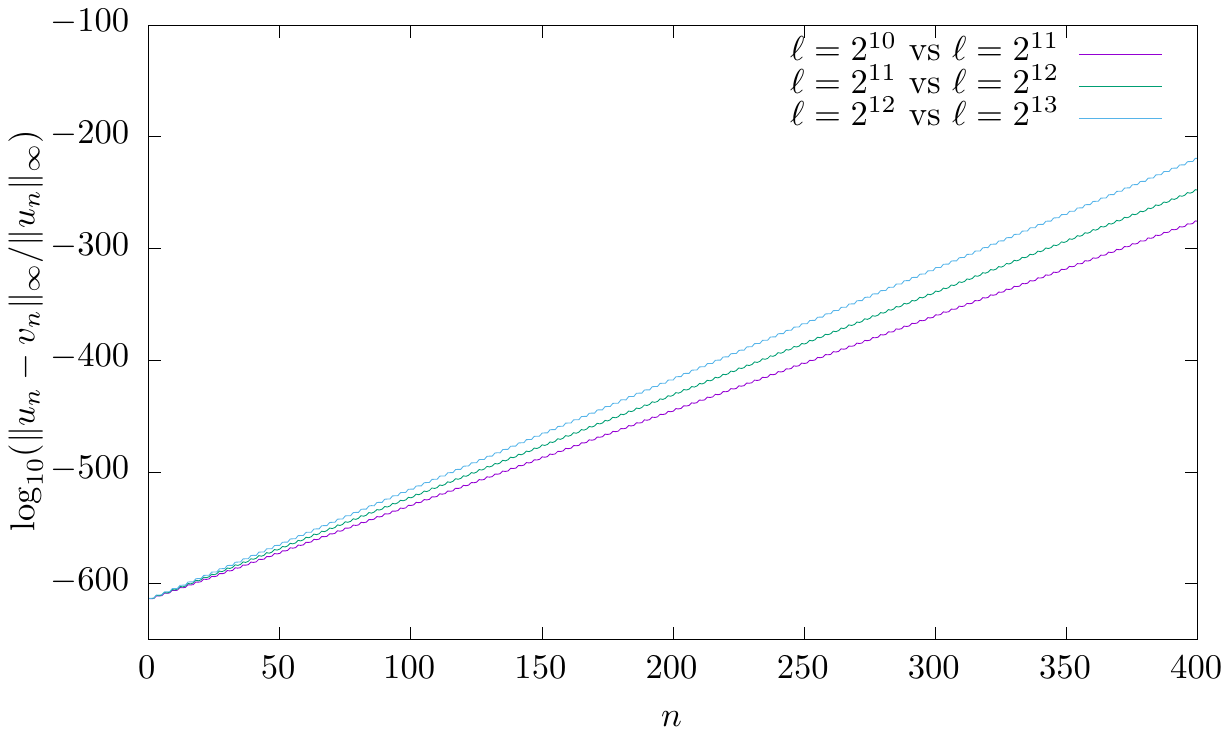} \hfil\hfil\hfil
\includegraphics[ trim = 8cm 10cm 8cm 10cm,  width=4truecm]{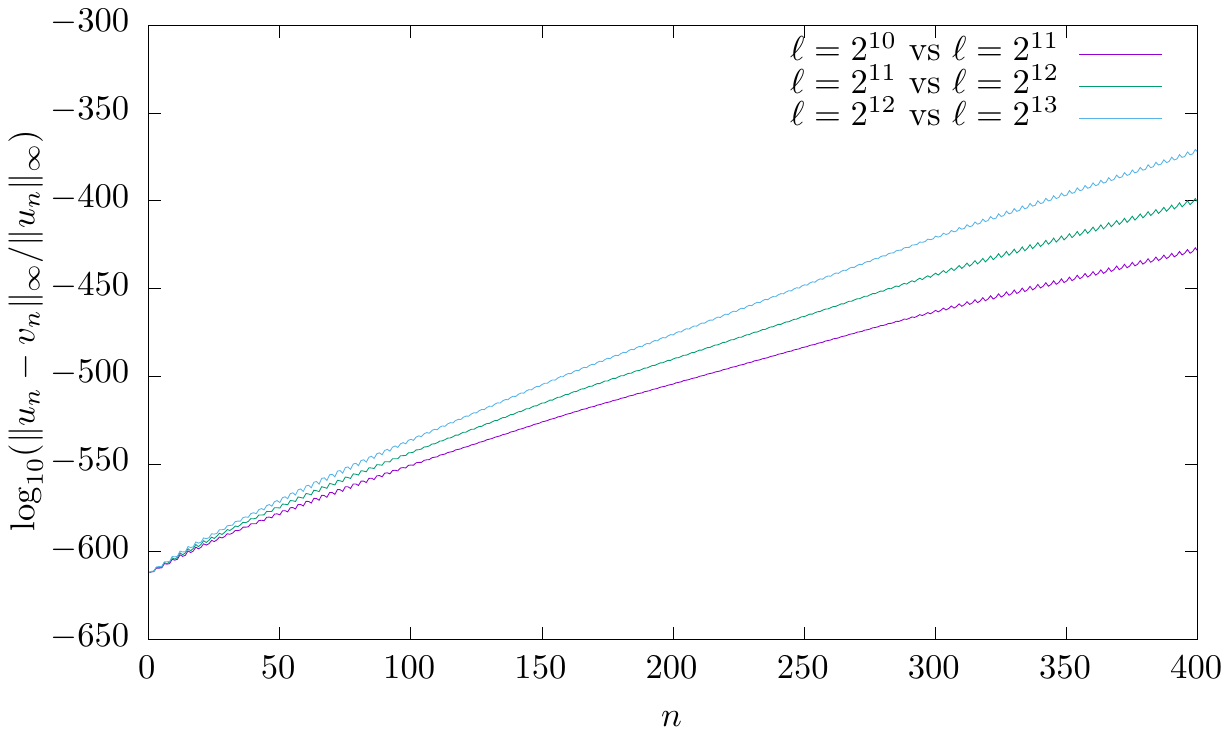}
\caption{Left panel: Graphs of $\log_{10}( \|u_n -v_n\|_\infty)$ where $u_n$ and $v_n$ are the same coefficients of the Lindstedt series but computed using a different number, $\ell$, of Fourier coefficients. Right panel: Graphs of the relative errors, $\log_{10}( \|u_n -v_n \|_\infty / \|u_n\|_\infty)$.} 

\label{errors}
\end{figure}

For this Corrigendum, the computations have been performed using 600 digits and  $2^\ell$ Fourier coefficients, with $10\leq \ell \leq 13$. 
Using this precision we have verified that the coefficients 
$u_n$ of the Lindstedt expansion have a relative error less 
than $10^{-300}$ when $n\leq 400$, see Figure \ref{errors}. We have also checked that the functions $u_n$ are trigonometric polynomials of degree $n$, as predicted in~\cite{BustamanteL20}, up to an error less than $10^{-200}$ within the same range of parameters. All the computations were done in pari/gp, \cite{pari}.

\section{Acknowledgements}
We thank Rafael de la Llave for fruitful discussions. We also thank the computer support at Georgia Tech. 
A.P.B. was partially supported by NSF grant DMS-1800241 . R.C. was partially supported by UNAM-DGAPA, PAPIIT project IN-101020.

\bibliography{correcbib}
\bibliographystyle{alpha}

\Addresses
\end{document}